\documentclass[12pt]{article}

\usepackage{amsmath}
\usepackage{amssymb}
\usepackage{amsthm}
\usepackage{latexsym}
\usepackage{color}
\usepackage{graphicx}
\usepackage{color}
\usepackage{soul}

\DeclareSymbolFont{calletters}{OMS}{cmsy}{m}{n}
\DeclareSymbolFontAlphabet{\mathcal}{calletters}

\def\be{\begin{eqnarray}}
\def\ee{\end{eqnarray}}

\def\b*{\begin{eqnarray*}}
\def\e*{\end{eqnarray*}}

\newtheorem{Theorem}{Theorem}[section]
\newtheorem{Definition}[Theorem]{Definition}
\newtheorem{Proposition}[Theorem]{Proposition}

\newtheorem{Lemma}[Theorem]{Lemma}

\newtheorem{Remark}[Theorem]{Remark}

\makeatletter \@addtoreset{equation}{section}

\usepackage{color}

\def \E{\mathbb{E}}

\def \L{\mathbb{L}}
\def \P{\mathbb{P}}

\def \R{\mathbb{R}}

\def\Dc{{\cal D}}

\def\Mc{{\cal M}}

\def\Pc{{\cal P}}

\def \0{\mathbf{0}}

\def\1{{\bf 1}}
\def \proof{{\noindent \bf Proof.\quad}}

\def \ep{\hbox{ }\hfill{ ${\cal t}$~\hspace{-5.1mm}~${\cal u}$}}

\usepackage{caption}

 \usepackage{color}

\title{Monotone Martingale Transport Plans\\
       and Skorohod Embedding}

\author{Mathias Beiglb\"ock\thanks{TU Vienna, mathias.beiglboeck@tuwien.ac.at. Supported through P26736 and Y782-N25.}
            \and
            Pierre Henry-Labord\`ere\thanks{Soci\'et\'e G\'en\'erale, Global Market Quantitative Research, pierre.henry-labordere@sgcib.com.}
            \and
            Nizar Touzi\thanks{Ecole Polytechnique Paris, Centre de Math\'ematiques Appliqu\'ees,
        nizar.touzi@polytechnique.edu. This work benefits from the financial support of the ERC Advanced Grant 321111, and the Chairs {\it Financial Risk} and {\it Finance and Sustainable Development}. }}

\date{\today}

\begin{document}

\maketitle

\begin{abstract}
We show that the left-monotone martingale coupling is optimal for any given performance function satisfying the martingale version of the Spence-Mirrlees condition, without assuming additional structural conditions on the marginals.
We also give a new interpretation of the left monotone coupling in terms of Skorokhod embedding which allows us to give a short proof of uniqueness. 
\end{abstract}

\vspace{0.9em}

{\small

}

\section{Introduction}

The roots of  optimal transport as a mathematical field go  back to Monge \cite{Mo81} and Kantorovich \cite{Ka42} who established its modern formulation.  More recently it experienced a steep development prompted by Brenier's  theorem \cite{Br91} and the milestone PhD-thesis of McCann \cite{Mc94}. The field is now famous for its striking applications in areas ranging from mathematical physics and PDE-theory to geometric and functional inequalities. We refer to \cite{Vi03, Vi09, AmGi13} for recent accounts of the theory.

Very recently there has also been  interest in optimal transport problems where the transport plan must satisfy additional martingale constraints. Such problems arise naturally in robust finance, but are also of independent mathematical interest, for example -- mirroring classical optimal transport -- they have important consequences for the study of martingale inequalities (see e.g.\ \cite{BoNu13,HeObSpTo12,ObSp14}).  Early papers to investigate such problems include \cite{HoNe12, BeHePe12, GaHeTo12, DoSo12, CaLaMa14}, and this topic is commonly referred to as martingale optimal transport. In mathematical finance, transport techniques complement the Skorokhod embedding approach to model-independent/robust finance (we refer to \cite{Ob04,Ho11} for an overview, see \cite{BeCoHu14} for a link with optimal transport theory). 

In view of the central role taken by Brenier's theorem  in  optimal transport, it is an intriguing question to obtain an analogous result also in the martingale setup. In this direction,  \cite{BeJu16, HeTo13} have proposed a martingale version of Brenier's monotone transport mapping.   

Under certain structural properties of the underlying probability distributions, it was established in \cite{HeTo13} that the monotone martingale transport plan solves the given variational problem for all coupling functions satisfying the \emph{Spence-Mirrlees} condition. In the present paper we show that this is in fact true for general continuous distributions $\mu, \nu$.

  A basic fundamental result is that the monotone martingale transport plan is \emph{unique}. However the original derivation (\cite{BeJu16}) is intricate and relies on subtle properties of potential functions and a  delicate approximation procedure. We provide a short conceptual proof of this result. 
This is based on a new connection with the classical Skorokhod embedding problem which may be of independent interest.


\section{Martingale optimal transport}
On the canonical space $\Omega:=\R\times\R$, we denote by $(X,Y)$ the canonical process, i.e. $X(x,y)=x$ and $Y(x,y)=y$ for all $(x,y)\in\Omega$. We also denote by $\Pc_\R$ and $\Pc_\Omega$ the collection of all probability measures on $\R$ and $\Omega$, respectively. For fixed $\mu,\nu\in\Pc_\R$ with finite first moments, our interest is on the following subsets of $\Pc_\Omega$:
 \be
 \Pc(\mu,\nu)
 &:=&
 \big\{\P\in\Pc_\Omega: X\sim_\P\mu, Y\sim_\P\nu\big\},
 \\
 \Mc(\mu,\nu)
 &:=&
 \big\{\P\in\Pc(\mu,\nu): \E^\P[Y|X]=X,~\mu-a.s. \big\}.
 \ee
The set $\Pc(\mu,\nu)$ is non-empty as it contains the product measure $\mu\otimes\nu$. By a classical result of Strassen \cite{St65}, we also know that $\Mc(\mu,\nu)$ is non-empty if and only if $\mu\preceq\nu$ in convex order, i.e.
 \be\label{convexorder}
 \mu(g)\le\nu(g)
 &\mbox{for all convex function}&
 g.
 \ee
Throughout we assume that $c:\Omega\longrightarrow\R$ is a measurable coupling function with $c\le a\oplus b$ for some $a\in\L^1(\mu)$ and $b\in\L^1(\nu)$. Here $a\oplus b(x,y):=a(x)+b(y)$ for all $(x,y)\in\Omega$. Then $\E^\P[c(X,Y)]$ is a well-defined scalar in $\R\cup\{-\infty\}$.
The martingale optimal transport problem, as introduced in \cite{BeHePe12} in the present discrete-time case and \cite{GaHeTo13}, is defined by:
 \be
 \mathbf{P}(c)
 \;:=\;
 \mathbf{P}
 &:=&
 \sup_{\P\in\Mc(\mu,\nu)} \E^\P[c(X,Y)].
 \ee
This problem is motivated by the problem of model-free superhedging in financial mathematics:
 \begin{equation}
 \mathbf{D} 
 :=
 \mathbf{D}(c)
 := 
 \inf_{(\varphi,\psi)\in\Dc} \big\{\mu(\varphi)+\nu(\psi)\big\},
 \end{equation}
where, denoting $h^\otimes(x,y):=h(x)(y-x)$ for all $(x,y)\in\Omega$,
 \begin{equation}
 \Dc:=
 \big\{(\varphi,\psi):~\varphi^+\!\in\L^1(\mu),~\psi^+\!\in\L^1(\nu),
          ~\mbox{and}~
          \varphi\oplus\psi+h^\otimes\ge c,~\mbox{for some}~h\in\L^\infty(\R)
 \big\}.
 \end{equation}
The following result was established in \cite{BeHePe12}.

\begin{Theorem}\label{thm:dualitymart}
Let $\mu\preceq\nu\in\Pc_\R$, and assume $c\in{\rm USC}(\Omega)$ with $c\le a\oplus b$ for some lower semicontinuous functions $a\in\L^1(\mu)$, $b\in\L^1(\nu)$. Then $\mathbf{P}=\mathbf{D}$, and $\mathbf{P}=\E^{\P^*}[c(X,Y)]$ for some $\P^*\in\Mc(\mu,\nu)$.
\end{Theorem}

We also recall the variational result of \cite{BeJu16} which provides a characterization of the optimal martingale measure $\P^*$.

\begin{Theorem}[Monotonicity Principle]\label{thm:variational}
Let $\mu\preceq\nu\in\Pc_\R$. 
Then, if $\P^*\in\Mc(\mu,\nu)$ is a solution of $\mathbf{P}$, there exists a support (i.e.\ a Borel set $\Gamma\subset\Omega$ with $\P^*[(X,Y)\in\Gamma]=1$) such that for all $\P^0\in\Pc(\mu^0,\nu^0)$ , $\mu^0,\nu^0\in\Pc_\R$ with finite support, and ${\rm Supp}(\P^0)\subset\Gamma$, we have
 \b*
 \E^{\P^0}\big[c(X,Y)\big]
 \ge
 \E^{\P}\big[c(X,Y)\big]
 &\mbox{for all}&
 \P\in\Pc(\mu_0,\nu_0)
 ~\mbox{with}~
 \E^{\P}[Y|X]=\E^{\P^0}[Y|X].
 \e*
\end{Theorem}
The fact that the absence of a duality gap (as in Theorem \ref{thm:dualitymart}) implies a variational result similar to the one given in Theorem \ref{thm:variational} is well known in the transport literature, see e.g.~Villani's book \cite[p88]{Vi03}. Extensions of  Theorem \eqref{thm:variational} have been provided in \cite{Za14, BeGr14, BeNuTo15}, see \cite{BeCoHu14} for a variant applicable to Skorokhod problem.
In particular Zaev~\cite{Za14} obtains (as  a special  case of his results) a version of Theorem \ref{thm:variational} with a simple proof under the assumption that the duality $\mathbf{P}=\mathbf{D}$ holds. For the convenience of the reader we report the argument\footnote{While Zaev's argument is more direct and intuitive (in our opinion), the approach of \cite{BeJu16, BeGr14} applies also in cases where the duality $\mathbf{P}=\mathbf{D}$ might fail (e.g.\ if $c$ is not upper bounded in the sense of Theorem \ref{thm:dualitymart}) and we state  Theorem \ref{thm:variational} in this slightly more general form.} from Zaev~\cite{Za14}. 

\smallskip

\proof
Pick  sequences of admissible dual functions $\phi_n, \psi_n, h_n, n\geq 1$ such that $\mu(\phi_n)+\nu(\psi_n)\to \mathbf{P}$ and fix ${\P^*}\in \Mc(\mu, \nu)$ such that $\E^\P[c(X,Y)]=\mathbf{P}.$ 
Since $\phi_n\oplus\psi_n+h_n^\otimes \geq c$ and
$$\E^{\P^*}[\phi_n(X)+\psi_n(Y)+h_n(X)(Y-X)] =\mu(\phi_n)+\nu(\psi_n) \to \mathbf{P}= \E^{\P^*}[c(X,Y)]$$ it follows that $\phi_n\oplus\psi_n+h_n^\otimes $ tends to $c$ in $\|.\|_{L^1(\P^*)}$. Passing to a subsequence which we denote again by $n$ we find that this convergence holds pointwise on a set   $\Gamma$ with ${\P^*}(\Gamma)=1$. 

Assume that $\P^0\in\Pc(\mu^0,\nu^0)$ , $\mu^0,\nu^0\in\Pc_\R$ with finite support, and ${\rm Supp}(\P^0)\subset\Gamma$, and that
$
\P\in\Pc(\mu^0,\nu^0)$ satisfies 
 $\E^{\P}[Y|X]=\E^{\P^0}[Y|X].$ Note that then $\E^{\P}[h(X)(Y-X)]=\E^{\P^0}[h(X)(Y-X)]$ holds for an arbitrary function $h$.  We thus obtain 
\begin{align*}
 \E^{\P^0}[c(X,Y)]   &=   \lim_n \E^{\P^0}[\phi_n(X)+\psi_n(Y)+h_n(X)(Y-X)] \\
 &=  \lim_n \E^{\P}[\phi_n(X)+\psi_n(Y)+h_n(X)(Y-X)] \geq \E^{\P}[c(X,Y)], 
\end{align*}
hence $\Gamma$ is as required.
\ep

\vspace{5mm}

\section{Monotone transport plans}
\label{sect:left monotone}
\setcounter{equation}{0}
\setcounter{Theorem}{0}

The following definition stems from \cite{BeJu16}.

\begin{Definition}\label{def:monotone}
We say that $\P\in\Mc(\mu,\nu)$ is left-monotone (resp. right-monotone) if there exists a support (i.e.\ a Borel set $\Gamma\subset\R\times\R$ with $\P[(X,Y)\in\Gamma]=1$) such that for all $(x,y_0), (x,y_1), (x',y')\in\Gamma$ with $x<x'$ (resp. $x>x'$), it must hold that
$y'\not\in(y_0,y_1)$.
\end{Definition}

The relevance of this notion is mainly due to the following extremality result which states that monotone martingale coupling measures are optimal for a class of martingale transport problems.

\begin{Definition}\label{MSM}
We say that a function $c:\R\times\R\longrightarrow\R$ satisfies the martingale Spence-Mirrlees condition if $c$ is measurable,  continuously differentiable with respect to $x$, and $c_x(x,.)$ is strictly convex on $\R$  for  $x\in\R$.
\end{Definition}

The martingale Spence-Mirrless condition was introduced in \cite[Remark 3.15]{HeTo13}, cf.\ also \cite[Definition 4.5]{HoKl12}. 

\begin{Theorem}\label{prop:optimal==>monotone}
Assume that the performance function $c:\R\times\R\longrightarrow\R$ satisfies the martingale Spence-Mirrlees condition. Then, any solution $\P^*$ of the martingale transport problem $\mathbf{P}$ is left-monotone.

In particular, for all $\mu\preceq\nu\in\Pc_\R$, there exists a left-monotone transference plan.
\end{Theorem}

\proof
The last existence result is a consequence of the existence of a maximizer $\P^*$ of the problem $\mathbf{P}(c)$ for some performance function $c$ satisfying the conditions of Theorem \ref{thm:dualitymart}.

Let $\P^*$ be a solution of the martingale transport problem $\mathbf{P}$, and suppose to the contrary that $\P^*$ is not left-monotone. Let $\Gamma$ be an arbitrary support of $\P^*$. By definition of the notion of left-monotonicity, we may find scalars $x<x'$ and $y_0<y'<y_1$ such that $(x,y_0),(x,y_1),(x',y')\in\Gamma$. We then introduce:
 \b*
 \P^0
 :=
 \frac12\big[\lambda\delta_{(x,y_0)}+(1-\lambda)\delta_{(x,y_1)}\big]
 +\frac12\delta_{(x',y')},
 \\
 \P
  :=
 \frac12 \delta_{(x,y')}
 +\frac12\big[\lambda\delta_{(x',y_0)}+(1-\lambda)\delta_{(x',y_1)}\big],
 \e*
where $\lambda\in(0,1)$ is defined by
 \b*
 \lambda y_0 + (1-\lambda) y_1
 &=&
 y'.
 \e*
Clearly, Supp$(\P^0)=\{(x,y^0),(x,y^1),(x',y')\}\subset\Gamma$. Moreover, $\int \P^0(dx,.)=\int\P(dx,.)=\mu^0(dx):=\frac12(\delta_{x}+\delta_{x'})(dx)$. Similarly, $\int \P^0(.,dy)=\int\P(.,dy)=\nu^0(dy):=\frac12(\lambda\delta_{y_0}+\delta_{y'}+(1-\lambda)\delta_{y_1})(dy)$. This shows that $\P,\P^0\in\Pc(\mu^0,\nu^0)$. We also directly compute that
 $$
 \E^{\P^0}[Y|X=x]=\E^{\P}[Y|X=x']=\lambda y_0+(1-\lambda)y_1
 ~\mbox{and}~
 \E^{\P^0}[Y|X=x']=\E^{\P}[Y|X=x]=y',
 $$
which in view of the definition of $\lambda$, shows that $\E^{\P^0}[Y|X]=\E^{\P}[Y|X]$. We may then apply the variational Theorem \ref{thm:variational}, and conclude that $\E^{\P^0}[c(X,Y)]\ge\E^{\P}[c(X,Y)]$, i.e.
 \b*
 g(x):=\lambda c(x,y_0)+(1-\lambda)c(x,y_1)-c(x,y')
 &\ge&
 g(x').
 \e*
We now show that this inequality is in contradiction with the martingale Spence-Mirrlees condition. Indeed, as $c$ is  continuously differentiable in $x$, we have
 \b*
 g(x')-g(x)
 &=& 
 \int_x^{x'} \big[\lambda c_x(\xi,y_0)+(1-\lambda)c_x(\xi,y_1)-c_x(\xi,y')\big] d\xi 
 \;>\;
 0,
 \e*
where the strict inequality follows from the strict convexity of the density $c_x$ in $y\in\R$.
\ep

\vspace{5mm}

For an atom-less measure $\mu$, the following easy consequence, reported from \cite{BeJu16}, shows that left (and right) monotone martingale transport plans have a very simple structure. Namely that the support is concentrated on two graphs. 

\begin{Proposition}\label{prop:leftmonotone2graphs}
Let $\mu\preceq\nu$ in convex order, and assume $\mu$ has no atoms. Let $\P\in\Mc(\mu,\nu)$ be a left-monotone transference plan. 
Then there exist functions $T_u, T_d:\R\longrightarrow \R$ with $T_d(x) \leq x \leq T_u(x)$, $\{T_d(x) =x\}=\{T_u(x)=x\}$, and such that 
 $$
 \P(dx,dy)
 =
 \mu(dx)\big[q(x)\delta_{T_u(x)}+(1-q(x))\delta_{T_d(x)}\big](dy)
 ~\mbox{with}~
 q(x):=\frac{x-T_d(x)}{T_u(x)-T_d(x)}\1_{\{T_d(x)<T_u(x)\}}.
 $$
Moreover, the pair $(T_d,T_u)$ is unique $\mu-a.e.$ and satisfies:
\begin{enumerate}
\item\label{Tu:increasing} $T_u$ is non decreasing.
\item\label{curtain} If $x<y$ then $T_d(y) \notin (T_d(x), T_u(x))$. 
\end{enumerate}
\end{Proposition} 

Given the structure of left monotone transport plans for an atom-less probability measure $\mu$, we call such a martingale transport measure a {\it left monotone transport map}.

Proposition \ref{prop:leftmonotone2graphs} should be compared to the situation in standard optimal transport problem, where the Fr\'echet-Hoeffding coupling defines a transport plan concentrated on a graph which solves simultaneously the Monge and the Kantorovitch problem. Clearly, one can not expect that martingale transport maps be concentrated on one single graph, as this would imply that the martingale is deterministic, and therefore constant which can happen only in the degenerate case $\mu=\nu$.

The following uniqueness result for the martingale transport problem $\mathbf{P}$ is a direct consequence of Proposition \ref{prop:leftmonotone2graphs}.

\begin{Proposition}[{cf.\ \cite[Section 5.1]{BeJu16}}]\label{prop:uniqueP}
Let $\mu,\nu\in\Pc_\R$ be such that $\mu\preceq\nu$ in the convex order, and $\mu$ without atoms. Let $c:\R\times\R\longrightarrow\R$ be a performance function satisfying the martingale Spence-Mirrlees condition. Then, there exists at most one solution to the martingale transport problem $\mathbf{P}(c)$.
\end{Proposition}

\proof
Suppose that $\mathbf{P}(c)$ has two solutions $\P_1$ and $\P_2$ in $\Mc(\mu,\nu)$. By Thorem \ref{prop:optimal==>monotone} and Proposition \ref{prop:leftmonotone2graphs}, it follows that both $\P_1$ and $\P_2$ are left monotone and concentrated on two graphs. We next consider the probability measure $\bar\P:=(\P_1+\P_2)/2$. Clearly, $\bar\P\in\Mc(\mu,\nu)$, and $\E^{\bar\P}[c(X,Y)]=\E^{\P_1}[c(X,Y)]=\E^{\P_2}[c(X,Y)]$ so that $\bar\P$ is also a solution of $\mathbf{P}(c)$. However, this is in contradiction with Theorem \ref{prop:optimal==>monotone} and  Proposition \ref{prop:leftmonotone2graphs}, as $\bar\P$ is concentrated on more than two graphs.
\ep

\section{Uniqueness  of the left-monotone transference map}

It was established in \cite{BeJu16} that for fixed marginals $\mu, \nu$ there exists a \emph{unique} left-monotone transference map in $\Mc(\mu,\nu)$. As the original argument is rather lengthy and maybe not entirely transparent, it seems worthwhile to revisit this basic (but important) result. 

In the rest of this paper, we assume that $\mu$ has no atoms and we fix $\P$ and the corresponding $T^\P:=(T_d, T_u)$ as in Proposition \ref{prop:leftmonotone2graphs}. Our aim is to prove that $\P$ is the unique left-monotone transference plan in $\Mc(\mu,\nu)$.

To this end, we use that a left-monotone martingale transport plan gives rise to a particular solution of the Skorokhod embedding problem:  for measures $\mu$, $\nu$ in convex order the Skorokhod problem is to construct a stopping time $\tau$ such for a  Brownian motion started in $B_0\sim \mu$, the distribution of $B_\tau$ equals $\nu$. We refer to \cite{Ob04, Ho11} for recent surveys on the Skorokhod embedding problem.

Typically one is interested to find a  minimal solution of the Skorokhod problem, i.e.\ a stopping time which in addition to $B_\tau\sim \nu$ satisfies that for all $\sigma \leq \tau$ with $B_\sigma\sim\nu$ one has $\tau=\sigma$. This notion was introduced by Monroe \cite{Mo72} in the case $\mu=\delta_0$, and further extended to the case of a general starting law by Cox \cite{Co08}. 

We first observe that the pair of maps $T^\P=(T_d,T_u)$ introduced in Proposition \ref{prop:leftmonotone2graphs} suggests to introduce the stopping time
 \be
 \tau^\P
 &:=&
 \inf\big\{ t>0:~B_t \not\in \big[ T_d(B_0),T_u(B_0) \big] \big\}.
 \ee

\begin{Proposition}\label{prop:min}
Let $B$ be a Brownian motion started from $B_0\sim\mu$. Then $(B_0,B_{\tau^\P}) \sim \P$, and $\tau^\P$ is a minimal stopping time.
\end{Proposition}

\proof
The fact that $(B_0,B_{\tau^\P}) \sim \P$ is clear by construction. Next, since $\tau^\P$ is a hitting time, we now show that it is a minimal embedding of $\nu$ starting from $\mu$. Minimality of hitting times with starting measure $\delta_0$ was observed by Monroe \cite{Mo72} (just after his Definition 1), we report the full argument here in order to emphasize that this result extends trivially to an arbitrary starting measure $\mu$.

Denote $A:=\big(T_d(B_0),T_u(B_0)\big)$. For a stopping time $\sigma\le\tau^\P$, we have $\{B_{\tau^\P}\notin A\}\supseteq\{B_\sigma\notin A\}$. Then, if $\sigma$ also embeds $\nu$, it follows that $\{B_{\tau^\P}\notin A\} = \{B_\sigma\notin A\}$, and equivalently $\{B_{\tau^\P}\in A\}=\{B_\sigma\in A\}$. Since $\tau^\P$ is the first exit time from $A$, this implies that, on the event set $\{\sigma<\tau^\P\}$, we have $B_{\tau^\P}\notin A$ while $B_\sigma\in A$. Then $\{B_{\tau^\P}\in A\}=\{B_\sigma\in A\}$ implies that $\tau^\P=\sigma$ a.s. 
\ep

\vspace{5mm}

In the context of the Skorokhod embedding problem 
a \emph{barrier} is  a measurable subset of $\R\times \R$ such that for any point $(x,y)$ contained in the barrier,  the whole line  $[x, \infty)\times \{y\}$ is a subset of the barrier.

We use the mappings $T_u, T_d$ to define the barrier
 \be\label{def:R}
 R:=R_u\cup R_d
 &\mbox{with}&
 R_i:= \bigcup_{x\in \R} [T_i(x)-x, \infty) \times \{T_i(x)\},
 ~~i\in\{d,u\}.
 \ee
 
\begin{center}
\begin{figure}[th]
\resizebox{.9\textwidth}{!}{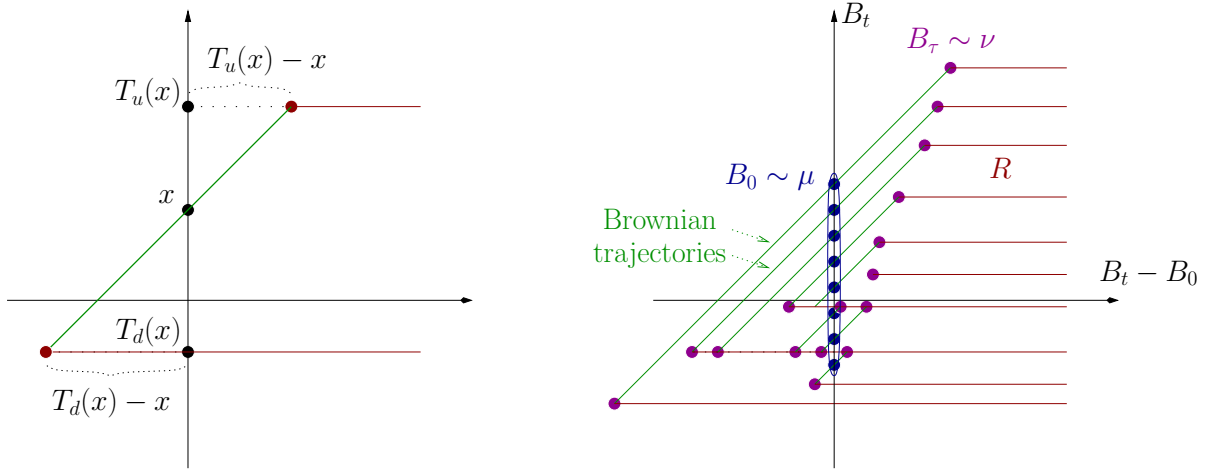}
 \caption{{\small\emph{ The left-monotone coupling as a barrier-type embedding:
 The left part depicts the construction of the barrier $R$. The right hand side shows how the set $R$ gives rise to a barrier-type stopping time in this particular phase-diagram. }}}
  \label{fig:RootSolution}
\end{figure}
\end{center}

We next introduce another stopping time defined as the first hitting time of this barrier:
 \be\label{LCtime}
 \tau_R
 &:=& 
 \inf \{t: (B_t-B_0, B_t)\in R\}.
\ee
It is clear that $\tau_R\le \tau^\P$. The following crucial result shows that equality in fact holds.

\begin{Lemma}\label{lem:tauR=tauP}
The stopping times $\tau_R$ and $\tau^\P$ are equal. In particular $\tau_R$ is a minimal stopping time and $(B_0,B_{\tau_R}) \sim \P$. 
\end{Lemma}

\proof
We only focus on the non-trivial inequality $\tau_R\ge \tau^\P$. To see this, we shall verify that the Brownian path started in $x$ hits $R$ either in $(T_u(x)-x, T_u(x))$ or in $(T_d(x)-x, T_d(x))$. Indeed, assume that a Brownian path hits a different line of the barrier $B_{\tau_R}\in [T_d(y)-y, \infty)\times\{T_d(y)\}$ or $B_{\tau_R}\in [T_u(y)-y, \infty)\times\{T_u(y)\}$, for some $y\neq x$. From our construction (see Figure \ref{fig:RootSolution}), we observe that, necessarily $y>x$, and

$\bullet$ in the case $B_{\tau_R}\in[T_d(y)-y, \infty)\times\{T_d(y)\}$, we have $T_d(x)<T_d(y)<T_u(x)$, contradicting Property \ref{curtain} of  Proposition \ref{prop:leftmonotone2graphs}.

$\bullet$ in the case $B_{\tau_R}\in [T_u(y)-y, \infty)\times\{T_u(y)\}$, we have $T_u(y)<T_u(x)$, contradicting Property \ref{Tu:increasing} of  Proposition \ref{prop:leftmonotone2graphs}.

Hence, $\tau_R=\tau^\P$, and the minimality property follows from Proposition \ref{prop:min}.
\ep

\vspace{5mm}

We have thus obtained an interpretation of the left-monotone transport plan in terms of a barrier-type solution to the Skorokhod problem.  This interpretation is useful for our purpose since  it allows us to use a short argument of Loynes \cite{Lo70} (which in turn builds on Root \cite{Ro66}) to show that there is only one left-monotone transference plan.

\begin{Lemma}[cf.~Loynes \cite{Lo70}]
Let $\P_1$, $\P_2$ be left-monotone martingale transport plans in $\Mc(\mu,\nu)$, with corresponding maps $T^\P=(T^i_u, T^i_d)$ satisfying the conditions of Proposition \ref{prop:leftmonotone2graphs}, and denote by $R^{\P_i}, i=1,2$ the corresponding barriers defined as in \eqref{LCtime}. Then $\tau_{R^{\P_1}}= \tau_{R^{\P_2}}$, a.s. 
\end{Lemma} 

\proof
For a set $A\subset\R$, we abbreviate $R_i(A):= R^{\P_i}\ \cap\ (\R\times A)$, $i=1,2$. Denote
 \b*
 K:= \big\{x: m_1(x) < m_2(x) \big\}   
 &\mbox{where}&
 m_i(x):=\inf\{m: (m,x) \in R^{\P_i}\},~~i=1,2.
 \e*
Fix a trajectory $(B_t)_t=(B_t(\omega))_t$ such that $B_{\tau_{R^{\P_2}}} \in K$. Then $(B_t-B_0,B_t)_t $ hits ${R_2}(K)$ before  it enters $R_2(K^C)$. 
But then $(B_t-B_0,B_t)_t  $ also hits $ R_1(K)$ before it enters $R_1(K^C)$. Hence 
$$
B_{\tau_{R^{\P_2}}} \in K \quad \Longrightarrow \quad B_{\tau_{R^{\P_1}}} \in K.
$$ 
As both stopping times embed the same measure, this implication  is an equivalence almost surely, and we may set $\Omega_K:=\{B_{\tau_{R^{\P_1}}} \in K\}= \{B_{\tau_{R^{\P_2}}} \in K\}$. On $\Omega_K$ we have $\tau_{R^{\P_1}} \leq \tau_{R^{\P_2}}$ while $\tau_{R^{\P_1}} \geq \tau_{R^{\P_2}}$ on $\Omega_K^C$. Then, for all Borel subset $A\subset\R$:
 \b*
 \P\big[B_{\tau_{R^{\P_1}} \wedge \tau_{R^{\P_2}}}\in A\big]
 &=&
 \P\big[B_{\tau_{R^{\P_1}}}\in A ,B_{\tau_{R^{\P_1}}}\in K\big]
 +\P\big[B_{\tau_{R^{\P_2}}}\in A ,B_{\tau_{R^{\P_1}}}\in K^c\big]
 \\
 &=&
 \P\big[B_{\tau_{R^{\P_1}}}\in A ,B_{\tau_{R^{\P_1}}}\in K\big]
 +\P\big[B_{\tau_{R^{\P_2}}}\in A,B_{\tau_{R^{\P_2}}}\in K^c\big]
 \\
 &=&
 \P\big[B_{\tau_{R^{\P_1}}}\in A\big],
 \e*
since $B_{\tau_{R^{\P_i}}}\sim\nu$. Hence $\tau_{R^{\P_1}} \wedge \tau_{R^{\P_2}}$ embeds $\nu$. Similarly, we see that  $\tau_{R^{\P_1}} \vee \tau_{R^{\P_2}}$ also embeds $\nu$. Since $\tau_{R^{\P_1}}$ and $\tau_{R^{\P_2}}$ are both minimal embeddings, we deduce that $\tau_{R^{\P_1}} \wedge \tau_{R^{\P_2}}=\tau_{R^{\P_1}} \vee \tau_{R^{\P_2}}$ and thus $\tau_{R^{\P_1}} = \tau_{R^{\P_2}}$. 
\ep

\vspace{5mm}

As $T_u, T_d$ can be $\mu$-a.s.\ recovered from the stopping time $\tau_R$ it follows that $T_u, T_d$ are uniquely determined and we obtain

\begin{Theorem}\label{Unique} Assume that $\mu \preceq \nu$ and that $\mu$ has no atoms. 
There exists precisely one left-monotone martingale coupling for  $\mu$, $\nu$.
\end{Theorem}

\begin{Remark} In \cite{BeJu16}, Theorem \ref{Unique} is proved without the assumption that $\mu$ has no atoms. It is possible to use the present approach to establish also this more general result; the basic idea is to represent the measure $\mu$ as an atom-less measure $\bar \mu$ on the set
\[ L:=\bigcup_{x\in {\rm Supp}\  \mu}\{x\}\times [0,\mu(\{x\})].\]
Since the extension of the above to arguments to the more general case is straightforward and the result is known from \cite{BeJu16} we do not elaborate. 
\end{Remark}

We conclude the paper with an additional property of our Skorohod embedding interpretation of the monotone martingale transport plan.

\begin{Proposition}
The process $B_{\wedge\tau_R}$ is a uniformly integrable martingale.
\end{Proposition} 

\proof
By Proposition \ref{prop:min} and Lemma \ref{lem:tauR=tauP}, $\tau_R=\tau^\P$ is a minimal embedding of $\nu$ with starting measure $\mu$.
Since $\E\big[B_{\tau^\P}\big]=\E[B_0]$, we know from Lemma~12 and Theorem 17 in Cox \cite{Co08} that minimality is equivalent to the uniform integrability of the process $B_{\wedge\tau^\P}$. 

For the convenience of the reader, we also provide a direct justification of the uniform integrability in our setting. Observe that $B_{\wedge\tau^\P}\in \big[ T_d(B_0),T_u(B_0) \big]$, a.s. Then, conditional on $B_0$, the process $B_{\wedge\tau^\P}$ is a bounded martingale. By the Jensen inequality, this provides
 \b*
 \E\big[\phi(B_{t\wedge\tau^\P})|B_0\big]
 \le 
 \E\big[\phi(B_{\tau^\P})|B_0\big],
 ~\mbox{a.s. for all convex function}~\phi.
 \e*
In particular, it follows that for any constant $c>0$:
 \b*
 \E\big[ |B_{t\wedge\tau^\P}|\1_{|B_{t\wedge\tau^\P}|\ge c}\big]
 &\le&
 \E\big[ \big(2|B_{t\wedge\tau^\P}|-c\big)^+\big]
 \;\le\;
 \E\big[ \big(2|B_{\tau^\P}|-c\big)^+\big],
 \e*
which provides the uniform integrability of the process $B_{\wedge\tau^\P}$. 
\ep
%
%
%

\end{document}